\documentclass[10pt,a4paper,twoside]{article}
\usepackage{graphicx}
\usepackage{amssymb}
\usepackage{amsmath}

\begin{document}

\title{The Mathematician's Bias -- and the Return to Embodied Computation}
\author{S. Barry Cooper
\\ {School of Mathematics, University of Leeds,}\\ {Leeds LS2 9JT, United Kingdom}
}
\date{}
\maketitle
\begin{abstract}There are growing uncertainties surrounding 
the classical model of computation established by 
G\"odel, Church, Kleene, Turing and others in the 1930s onwards. 
The mismatch between the Turing machine conception,  and the experiences 
of those more practically engaged in computing, has parallels with 
the wider one between science and those working creatively or intuitively 
out in the `real' world. The scientific outlook is more flexible and 
basic than some understand or want to admit. The science is subject to 
limitations which threaten careers. We look at embodiment and 
disembodiment of computation as the key to the mismatch, and 
find Turing had the right idea all along -- amongst a  
productive confusion of ideas about computation in the real 
and the abstract worlds. 
\end{abstract}

When we get out of bed in the morning, we approach a complicated  
world of information with a determination not just to survive the day -- 
though that may be hard enough: we mean to ``compute" our way towards 
various vaguely defined objectives. The process will likely be messy, but 
we certainly experience it as a computational one. Our conception 
of the computation is very flexible, but we host no {\em in principle} 
rejection of Turing's notion of mechanical intelligence.  

But our sense of ownership of a computational process deserts 
us somewhat when we think about what it is that makes our 
daily computing adventure so 
complicated. The world outside has both predictability and a lack of it 
bordering on randomness. 
Here is Nassim Taleb,  in his best-selling book ``The Black Swan" \cite{Ta07}: 
\begin{quote}
I have spent my entire life studying randomness, practicing 
randomness, hating randomness. The more that time passes, the worse things 
seem to me, the more scared I get, the more disgusted I am with Mother 
Nature. The more I think about my subject, the more I see evidence that the 
world we have in our minds is different from the one playing outside. 
Every morning the world appears to me more random than it did the day 
before, and humans seem to be even more fooled by it than they were 
the previous day. It is becoming unbearable. I find writing these lines 
painful; I find the world revolting.
\end{quote}

Taleb distrusts mathematicians and their models, from personal experience of 
their failures, and of their perceived unwillingness to face up to the realities. 
But we need to give the professionals a chance. Let us look more closely 
at the modelling process, and how we deal with computing 
in a material world. 
Can we absorb Taleb's computational context into a classical model based 
on logical structure. Or does embodied information need to be separately 
modelled? And does this take us beyond the mathematician's 
focus on computable {\em functions}.

\section{Computation Disembodied}

What was clearly new about the Turing model of computation 
was its successful {\it disembodiment} of the machine. For 
practical purposes, this was not as 
complete as some post-Turing theoreticians like to pretend: the re-embodied computer 
which is now a familiar feature of the modern world was hard won 
by pioneering engineers. But, for the purposes of the stored 
program computer, and for the proof of incomputability of 
the Halting Problem, the essential disembodiment was that delivering  
program-data convergence. It was this universality that 
John von Neumann recognised as a theoretical anticipation of the stored 
program computer. The apparent omnipotence even led Turing 
to talk of his post-war ACE project as aimed at building `a brain'.

This paradigm has achieved a strong grip on subsequent thinking. 
Within the philosophy of mind there is a strong tendency towards physicalism and 
functionalism, both of which open the door to some version of the 
Turing model. The functionalist (see Hilary Putnam \cite{Put60}) stresses what a computer 
{\it does} as something realisable in different  hardware. 
 An important expression of the functionalist view in computer 
science is provided by the notion of a {\it virtual machine}, whereby 
one expects to achieve software implementation of a 
given programmable machine. Aaron Sloman \cite{Slo09} and others have usefully  
applied the concept to AI.

This playing down of distinction between 
information and process has been taken further, and become a familiar feature of 
programming and theory. As Samson Abramsky describes  
(private communication):
\begin{quote}
Turing took traditional mathematical objects, real numbers, functions 
etc. as the things to be computed. In subsequent work in Computer Science, 
the view of computation has broadened enormously. In the work on 
concurrent processes, the behaviour is the object of interest. There is 
indeed a lack of a clear-cut Church-Turing thesis in this wider sphere of 
computation -- computation as interaction, as Robin Milner put it.
\end{quote}
In the quantum world there is a parallel 
convergence  between matter and law-like energy. All this has given rise to 
a standard computational paradigm   vulnerable to 
surprises from the natural world. Physical processes 
not subject to data-process convergence will not be {\it recognisably} 
different. But beneath the `normal science', 
theoretical inadequacies may be brewing -- or not, 
according to the viewpoint. The challenges to the standard model 
are varied, but most seem to have an impact on universality. 

What is happening in both situations is a side-stepping of the {\em mathematically} 
familiar type structure, whereby numbers, functions and relations, and relations over functions etc. 
give rise to a hierarchy of fundamentally different objects, increasingly hard to handle 
as one closes up hierarchically the universe of definable entities. In the real world 
we require a level of constructibility, of computability even, which forces approximations 
and a rejection of paths to the unknown.

\subsection{The Mathematician's Bias}

A symptom of the inadequacy of a type-constrained world-view is the 
October 2010 ACM Ubiquity Symposium on {\it What is Computation?} 
Part of the Editor's Introduction by Peter J. Denning \cite{Denning:2010} reads:
\begin{quote}
By the late 1940s, the answer was that computation was steps 
carried out by automated computers to produce definite outputs. 
That definition did very well: it remained the standard for nearly fifty 
years. But it is now being challenged. People in many fields have accepted 
that computational thinking is a way of approaching science and engineering. 
The Internet is full of servers that provide nonstop computation endlessly. 
Researchers in biology and physics have claimed the discovery of natural 
computational processes that have nothing to do with computers. How 
must our definition evolve to answer the challenges of brains computing, 
algorithms never terminating by design, computation as a natural occurrence, 
and computation without computers?
\end{quote}

In another contribution to the Symposium, Lance Fortnow \cite{Fortnow:2010} asks: ``So 
why are we having a series now asking a question that was settled in the 1930s?"
And  continues:
 \begin{quote}
 A few computer scientists nevertheless try to argue that the 
 [Church-Turing] thesis fails to capture some aspects of computation. Some 
 of these have been published in prestigious venues such as Science, the 
 Communications of the ACM  and now as a whole series of papers in ACM 
 Ubiquity. Some people outside of computer science might think that there is a 
 serious debate about the nature of computation. There isn't.
 \end{quote}

 Undeterred, Dennis J. Frailey thinks it's the mathematicians have got it wrong:
  \begin{quote}
 The concept of computation is arguably the most dramatic
advance in mathematical thinking of the past century \dots 
Church, G\"odel, and Turing defined it in terms of mathematical functions \dots  \quad
 They were inclined to the view that only the algorithmic functions constituted computation. 
I'll call this the ``mathematician's bias" because I believe it limits our thinking and 
 prevent us from fully appreciating the power of computation. 
 \end{quote}
 
Clearly, we do not have much of a grip on the issue. 
It is the old story of the {\it Blind Men and the Elephant} again. 
On the one hand computation is seen as an essentially open, 
contextual, activity with the nature of data in question. Others 
bring out a formidable armoury of mathematical weapons 
in service of the reductive project -- for them there is little  
in concurrency or interaction or continuous data or mental-recursions to 
stretch the mathematical capabilities of the Turing machine. 

Of course, there  has always been creative play on the paradoxical 
misfit between the Turing machine model and the realities 
of the real world.  One of the best-known is Jin Wicked's 
wonderful image of Alan Turing himself embodying his universal 
machine. How well this fits with David Leavitt's  
perception, in his book {\em The Man Who 
Knew Too Much},  of Turing actually {\em identifying} with 
his computing machines: 

\begin{figure}[!ht] \label{fig}
\begin{center}
\includegraphics[scale=0.28]{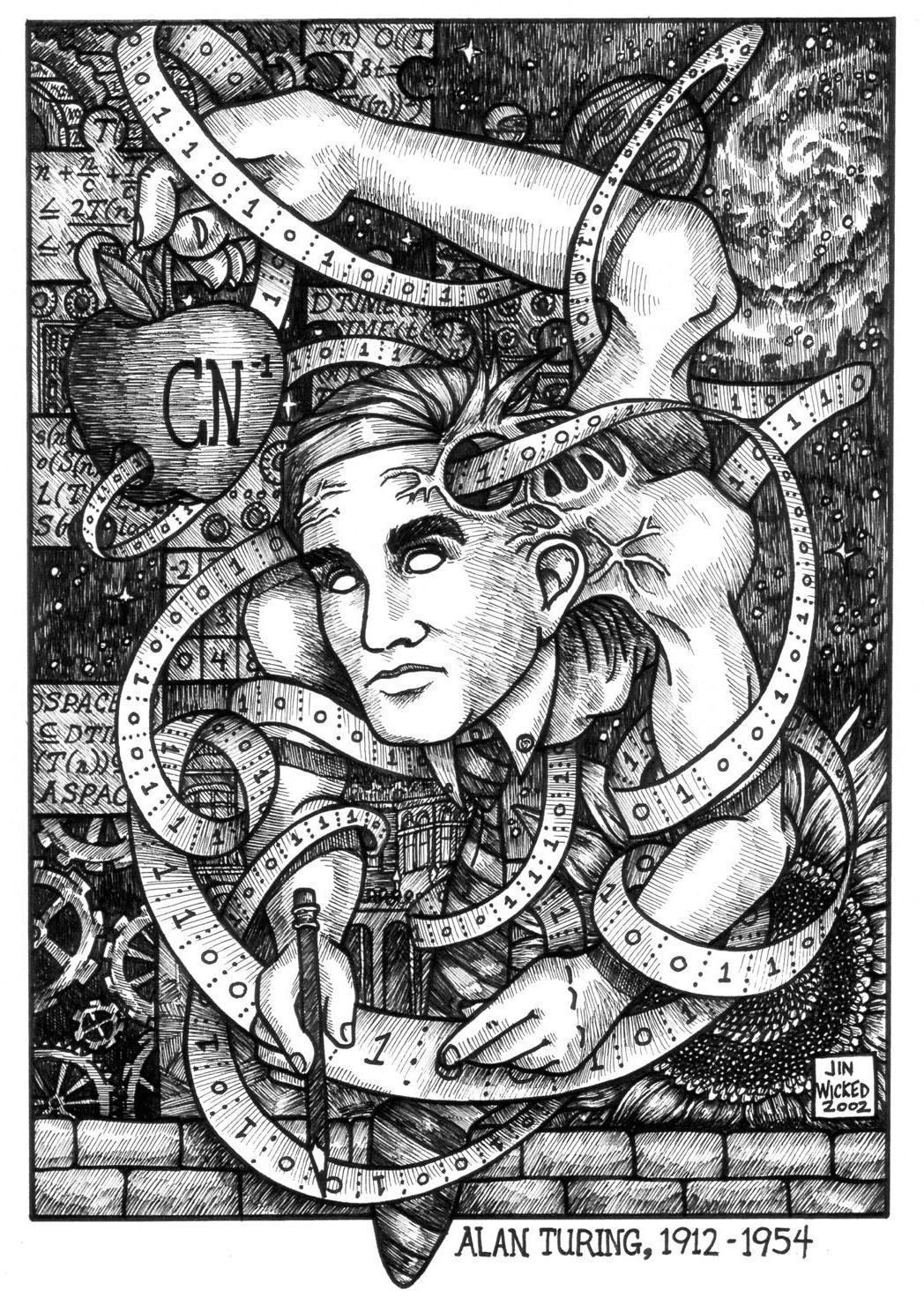}
\end{center}
\caption{Alan Turing as a Universal Turing Machine.\label{fig1}}
\end{figure}

Essentially, what is happening is that some observers are reviving the mathematical 
type-structure in a real-world context, others are denying that the mathematics 
is capable of inhabiting the material world. Some are struck by the sheer 
globality of how the world computes, the computational fruits of 
complexity, connectivity and interaction. And by the loss of a simple 
inductive structure implicit in non-linearity and the failure 
of computable approximations 
implicit in too inclusive a view 
 of computationally based environments. 

In the past though, it has been the mathematics that 
has clarified difficult problems and vague intuitions. 
Apparently, there is a lot of mathematics we have not found out 
how to use in the material context, though the ownership of the mathematics 
of reality is beginning to slip from the grasp of  those 
unwilling to adapt. 

At the cutting edge there is new mathematics being developed, and questioning 
of previously sacrosanct conceptual frameworks. Here are some uncomfortable 
observations (private communication) from Samson Abramsky:
\begin{quote}
Formally, giving a program + data
logically implies the output (leaving aside non-determinism or 
randomness), so why actually bother computing the result! \dots\\
\dots Can information increase in computation? Information theory 
and thermodynamics seem to tell us that it can't, yet intuitively, this is 
surely exactly why we compute -- to get information we didn't have before.
\end{quote}
And in mathematics, our operations and definitions certainly do 
give us new information. Is our wonderful universe so 
constrained it cannot go where even  high-school 
arithmetics leads us (as the Davis-Matiyasevich-Putnam-Robinson 
negative solution to Hilbert's Tenth Problem tells us)?

\section{The Mathematics of Embodiment?}

So it is not that mathematicians are only interested in simple mathematical 
objects. Or that computability only deals with functions over the natural numbers. 
What we do have is a computational paradigm which dominates our view 
of the landscape. We even have higher-type computation in various forms, 
including that mapped out by Stephen Kleene in his three late papers on the topic. 
What we did not get was any suggestion that higher-type computability might 
play a role in modelling the {\em processes} which it is now suggested 
might stretch the old Turing model. And the suspicion is that this simple connection 
can only be explained in the context of a powerful counter-paradigm. 

The other by-product of the counter-paradigm is a tendency to desert basic 
mathematical theory in favour of less focused descriptive arguments for 
new computational phenomena, accompanied by attempts at models 
derived from these descriptions unrooted in any classical analysis 
of their power. Here is another of the {\it Ubiquity} symposium 
contributors, Peter Wegner 
writing with Dina Goldin \cite{WG2005} on {\it The Church-Turing Thesis: 
Breaking the Myth}: 
 \begin{quote}
One example of a problem that is not algorithmic is the following instruction 
from a recipe [quote Knuth, 1968]:\\ 
`toss lightly until the mixture is crumbly.'\\
This problem is not algorithmic because it is impossible for a computer to know how long to mix: this may depend on conditions such as humidity that cannot be predicted with certainty ahead of time. In the function-based mathematical worldview, all inputs must be specified at the start of the computation, preventing the kind of feedback that would be necessary to determine when it's time to stop mixing.
 \end{quote}
At the level of a human trying to carry out the cooking instructions, 
we have a strong impression that something non-algorithmic 
is going on. But a determined reductionist would 
not be at all convinced that this is more than a superficial 
testing of classical modelling based on a particular mode 
of observation. The recipe is certainly implementable. 
But will be executed differently by different cooks under different 
conditions. 
What is actually happening in a global sense can be made 
precise and potentially reducible to the Turing model 
by modelling the {\em whole} context, cook and 
kitchen and ingredients. Who cares that the language used in the 
description of the recipe 
is a bit imprecise?

Even quantum physics has received effective attention from the reductionists, via    
David Deutsch's \cite{Deutsch85} placing of the standard model of quantum computing 
firmly within the scope of the Turing model. Of course, the model does not make 
full use of wave function collapse, and Calude and 
Svozil \cite{CS08} have shown, under reasonable assumptions, that 
quantum randomness entails incomputability. Meanwhile, the mind is very hard to 
pin down theoretically, but Deutsch is following a well-established  reductive 
tradition when he argues (in {\it Question and Answers with David Deutsch}, 
on the New.Scientist.com News Service, December, 2006):
\begin{quote}
I am sure we will have [conscious computers], I expect they will be purely 
classical, and I expect that it will be a long time in the future. 
Significant advances in our philosophical understanding of 
what consciousness is, will be needed.
\end{quote}

Maybe starting at the other end, and attempting 
to embody the richness of trans-computational mathematics, can 
bring more convincing results? Here we have a very 
different problem. Our examples of incomputable 
objects are very simple, even `natural' from some 
perspectives. But there is a huge gap between the 
universal Turing machine, which we can see embodied 
in a very persuasive sense by a modern computer, and the 
halting problem, which is incomputable but very abstract 
to the point of having no visible embodiment at all. 
And the lazy perception is that such examples of incomputability 
are nothing to do with the material world, even if we are 
moved to suspect trans-Turing computation by our 
impressionistic observation 
of natural phenomena. 

The key to bringing some clarity to the situation 
is to examine the mathematics of the incomputable, 
in all its basic simplicity. And to look for 
qualitatively similar mathematics bringing with it a level 
of embodiment, and an apparent avoidance of 
computability, to give us a better handle on the 
natural world and its candidate 
trans-Turing phenomena. 

A meaningful first example, residing at the border between mathematics 
and the embodied world is the Mandelbrot set. In fact, there are many 
other such fractal-like objects. But the Mandelbrot set has attracted special 
attention for good reasons: It has a simple definition over the complex plane 
involving basic arithmetical operations and a couple of quantifiers -- in fact, 
with a little fine-tuning, just one universal quantifier; it is graphically beautiful 
and complex, being both approximable via a computer screen, and 
containing endlessly explorable inner structure; and, its 
computability as a set of complex numbers provides us with 
a challenging open problem. As Roger Penrose says in his 
best-selling 1994 book \cite{Penrose2002} 
{\em The Emperor's New Mind}: 
\begin{quote}
Now we witnessed \dots 
a certain extraordinarily complicated looking set, namely the Mandelbrot set. 
Although the rules which provide its definition are surprisingly simple, the set itself exhibits
an endless variety of highly elaborate structures.
\end{quote}
What we see is a complicated and undeniably physical 
object existing in front of us on our computer screen {\em because} 
it was described, and having complex form {\em because} 
the description had a slightly elevated level of complexity 
of language -- namely the addition of a quantifier. Mathematically, 
the description is not much more complicated than that of the 
non-halting set of a universal Turing machine. The Mandelbrot 
set is of course made up of complex numbers rather than natural 
numbers, but the fact it appears on out computer screen gives 
it a digital kinship with familiar Turing-machine data. 
Leaving this abstract unembodied vista behind, and walking round to 
the other side of the Mandelbrot set, we see why such fractals 
have such significance for so many people -- the scene is one of 
Nature in all its fascinating complexity, beautiful and embodied, like the 
Mandelbrot set, but this time produced by simple laws which are themselves 
embodied. This toy fractal provides a neat connection between natural 
complexity and the unrestrained type-structures of mathematics, 
with their power to take us definably far beyond what we can 
computably capture. 

What indicates to the observer a phenomenon of 
computational complexity? It is the richness of {\em visible form} 
which impresses. It is the character of a higher order entity. 
On the other hand, despite the appearance pointing to 
the basic computational unpredictability, there is an overall 
identity to the appearance which we are capable of 
appreciating, and which is caught by the {\em definition} 
of the exotic shape. And it is natural to view this 
definition as {\em computed} by the totality of the 
underlying computational context. The mathematics 
both points to the atomic unpredictability 
of the universal Turing machine observed; and 
reassures us with a computation of a higher order, 
whereby we reaffirm a level of understanding of 
the underlying turbulence. 

There is an obvious 
role here for some of the most abstract and little known 
mathematics logicians have developed. How many 
mathematicians have imagined any role at all 
 in the real world for the  
conceptual development of higher order 
computation by Stephen Kleene, Gerald Sacks, 
Dag Normann and their successors? 
As a primary schoolboy, I was fascinated by the family 
folding mangle, and carried it off to possess and enjoy 
in my own little `camp' in the undergrowth beyond the back 
gardens. The usefulness had no meaning, the loss I was made aware 
of later by the grown-ups was a shock. It is surely time for the 
generalised `recursion theorists' to return this 
conceptually beautiful work to the real world which 
indirectly gave rise to it. And to make sense of the sort of 
physical mysteries that Alan Turing himself identified 
as important sixty years ago.

\section{ Emergent Natural Patterns}

Back in the early 1950s, Alan Turing became interested 
in how certain patterns in Nature arise. His seminal 
importation \cite{Tu52} of mathematical techniques into this 
area were to bring a new level of computability to 
natural phenomena, while improving our understanding 
of some of the mysteries pointed to by those who distrust 
the reductionism of those trapped by the Turing machine 
paradigm. Turing was able to relate 
such familiar features of everyday life as patterns 
on animals' coats -- stripes on zebras, patches on cow hides, 
moving patterns on tropical fish -- to simple reaction-diffusion 
systems describing chemical reactions and diffusion. 
This brought mathematics to play in biology in a way 
that made his one published paper on the topic 
one of the most cited in the literature, and the most 
cited in recent years of all the papers  he wrote -- 
including the 1936 computable numbers paper, and the 
influential AI `Turing Test' article \cite{Turing1950} in Mind in 1950. 

Turing's mathematics 
links up with the relationship between emergence and 
descriptions pointed to by the fractal example. Of course, 
some halting problems are solvable, as  
some Julia sets are computable (see  \cite{BBY07}). 
See \cite{Co09a} for a more detailed argument for the 
two-way correspondence between emergence and appropriate  
descriptions, opening the door to incomputability 
in Nature. And for the grounding of the   observed
 robustness of  emergent phenomena in  
mathematical framing of the descriptions as  
mathematical {\it definability}. 

Turing's approach is seminal, illuminating the 
connection between possible incomputability, mathematics, 
and natural phenomena. 
It has been carried forward by 
James D. Murray \cite{Murray02} and  
others, and though  
things get a lot more  
complex than the examples tackled by Turing,  
it is enough to make something more coherent from 
the confusion of intuitions and models. 
The embodiment does 
extend to the emergent halting set and possibly hierarchically beyond, 
taking us into a world beyond basic algorithms -- 
see  Chaitin's  recent take \cite{Chaitinip} on creativity 
and biology. 

This knack Turing had for drawing  fundamental computational aspects  
out of  concrete contexts has become increasingly clear to us 
since he died. The relevance of this seminal work on 
morphogenesis becomes daily wider and deeper in import. 
And we can now see it as definability {\em embodied}, and 
hence identify definability as a  very real form of {\em computability}. 

To summarise: It is often possible to {\em define} emergent 
properties in terms of the elementary actions underlying them. 
While in mathematics, relations and even 
objects {\em arise} from descriptions via the notion of 
definability. And, if the language used is complicated enough, 
this can be a source of Turing  incomputability. The key 
observation is that all that fancy higher type mathematics  
which the logicians carted off into the bushes has a vitally  
important practical usefulness. Emergent phenomena not 
only generate descriptions -- there is no serendipity here -- 
but {\em derive} from them.

One should not be misled by the morphogenesis into 
thinking of definability as exclusively relating higher order 
relations to basic algorithmic structure. As is well-known from the 
mathematics, definability can bring other quite basic aspects 
of a structure 
into play via descriptions in terms of objects 
of similar complexity. At the most basic level, we may 
find individual objects attaining their identity -- that is, a 
unique role as observable entities -- via their context 
within a wider causal framework. 

The picture is one of simple computable rules, with a degree of connectivity 
underpinning a higher order computation, and emergent forms 
{\em defined} at the edge of computability. 
We see morphogenesis inhabit the same world a the Mandelbrot 
set; the same world as the halting problem for a Turing machine; 
occupy  the same world as the large scale structure 
we see in the wider 
universe; and, one speculates, the same world as human mental 
creativity, linking the synaptic connectivity of the brain via 
neural net modelling to the most surprising 
artistic and scientific achievements. 

\section{The Mind as Mathematics?}

There is a strong inclination amongst both computer scientists and 
philosophers towards some kind of physicalist basis  
for mentality. For the philosopher, this is still an area 
with little to agree upon -- except the {\em supervenience} 
of the mental upon the physicality of the brain. As Jaegwon Kim 
puts it in his book \cite{Kim2000} on {\em Mind in a Physical World} (pp.14--15), supervenience:
\begin{quote}
represents the idea that mentality is at bottom physically based, and that there is no free-floating mentality unanchored in the physical nature of objects and events in which it is manifested.
\end{quote}
Or, more mathematically, from the {\em Stanford Encylopedia of Philosophy}:
\begin{quote}
A set of properties $A$ supervenes upon another set $B$  
just in case no two things can differ with respect to $A$-properties without also differing
with respect to their $B$-properties.
\end{quote}
In this context, one has familiar questions, such as: How can mentality have 
a causal role in a world that is fundamentally physical? And 
the puzzle of causal `overdetermination' -- the problem 
of phenomena having both mental and physical causes. As Kim \cite{Kim2005} sums up 
the problem (in {\em Physicalism, or Something Near Enough}, 2005):
\begin{quote}
\dots the problem of mental causation is solvable only if mentality 
is physically reducible; however, phenomenal consciousness 
resists physical reduction, putting its causal efficacy in peril.
\end{quote}

For computer scientists -- and mathematicians such 
as Alan Turing -- a physicalist approach takes one in the direction of computational 
models based on what we know of brain functionality. 

The computational content of most of these models is not so 
radically different from that of a Turing machine, though the level of connectivity 
and parallelism present may present some problems for the reductionists.  
While back in the real world, the typical observer of brain 
functionality will probably not be so 
impressed, and will 
see a lot of what one knows about brains,  from having one, missing. 
There are various models relevant here, most of them with a lot 
in common, and discussions which are transferable from one to 
another.     Rodney Brooks 
\cite{Brooks01} comments:
\begin{quote}
\dots neither AI nor Alife has produced artifacts that could be confused with a 
living organism for more than an instant. AI just does not seem as present 
or aware as even a simple animal and Alife cannot match the complexities 
of the simplest forms of life.
\end{quote}
The mind as a computational instrument 
presents  a level of challenge to the modellers  only 
matched by that facing those wanting 
to put the standard model of particle physics. 
 on a more secure footing. And many of us are persuaded of 
 important connections between these modelling 
 challenges, explored in various ways by different researchers -- 
 with, for example, the contrasting proposals of Roger Penrose and Henry Stapp 
 having attracted  widespread attention. 
 
It is still true that we have very interesting connectionist models for the brain 
(see \cite{Teuscher02} for Turing's contribution), 
and we are not surprised to see a leading researcher like 
Paul Smolensky \cite{Smolensky88} saying:
\begin{quote}
There is a reasonable chance that connectionist 
models will lead to the development of new 
somewhat-general-purpose self-programming, 
massively parallel analog computers, and 
a new theory of analog parallel computation: 
they may possibly even challenge the 
strong construal of Church's Thesis as the 
claim that the class of well-defined computations 
is exhausted by those of Turing machines.
\end{quote}
But,  as  Steven Pinker puts it:  
 ``\dots neural networks 
alone cannot do the job" --  and we too  expect a bit more than the emergence of 
embodied incomputability,  occupying a different level of data to that 
used for further computation. That is no real advance on the familiar 
Turing machine. It is not very promising for trans-Turing computing 
machines. Pinker does not talk about incomputability, but 
does describe  \cite[p.124]{pinker} human thinking as exhibiting  ``a kind of mental 
fecundity called recursion", giving us a good impression 
of real life emergent phenomena re-entering the 
system -- or the computational process as we would see it. 

Is there really something different happening here? If so, 
how do we model it? Does this finally sink the 
standard Turing machine model? Neuroscience gives 
us an impressively detailed picture of brain functionality. 
Antonio Damasio \cite{Damasio1999} vividly fills out the picture we get from Pinker:
\begin{quote}
As the brain forms images of an object - such as a face, a melody, a 
toothache, the memory of an event - and as the images of the object 
affect the state of the organism, yet another level of brain 
structure creates a swift nonverbal account of the events that are 
taking place in the varied brain regions activated as a 
consequence of the object-organism interaction. The mapping 
of the object- related consequences occurs in first-order neural 
maps representing the proto-self and object; the account of the 
causal relationship between object and organism can only be captured in second-order neural maps. 
\dots one might say that the swift, second-order nonverbal account 
narrates a story: that of the organism caught in the act of representing 
its own changing state as it goes about representing something else.
\end{quote}
The book gives a modern picture of how the human body 
distributes its `second-order' representations across the 
impressive connectivity of the human organism, and 
enables representations to play a role in further thought 
processes. 

Remarkably, once again, we find clues to what is happening mathematically in 
the 1930s work of Kleene and Turing. In his 1939 Princeton 
paper, Turing \cite{Turing1939} is in no doubt he is saying something 
mathematically about human mentality:
\begin{quote}
Mathematical reasoning may be regarded \dots as the exercise of a
combination of \dots intuition and ingenuity. \dots In pre-G\"odel times it 
was thought by some that all the intuitive judgements of mathematics 
could be replaced by a finite number of \dots rules. The necessity for 
intuition would then be entirely eliminated. In our discussions, 
however, we have gone to the opposite extreme and eliminated 
not intuition but ingenuity, and this in spite of the fact that our aim 
has been in much the same direction.
\end{quote}

To cut things short, Turing's analysis depends on Kleene's trick  
of building the constructive ordinals by computably {\it sampling} 
infinitary data needed to compute higher ordinals. The sampling 
is not canonical, there are lots of ways of doing it computably. 
This means we have an iterative process analogous to that 
described by Pinker and Damasio, which exhibits the 
irreversibility we expect from Prigogine. It is a picture of 
definable data represented using constructive ordinals. 
Turing's representation of data is scientifically standard, 
in terms of reals. 

Buried away in this opaque but wonderful paper is a key 
idea, that of the oracle Turing machine, which remarkably 
is just what one needs to model computable 
`causality' in science. One can represent any of the 
familiar physical transitions capturable on our 
computers -- and this encompasses a comprehensive 
swathe of what we regard as scientific -- using the appropriate kind 
of {\em relative} computation executable by the right oracle machine. 
The basic computational structure of 
the connectivity is captured by {\it functionals} modelled 
by oracle Turing machines. The mathematics delivers a complex 
structure -- an extended Turing model -- with a rich overlay of definable relations, 
corresponding to real world ubiquity of emergent form 
impacting non-trivially on the development of the 
organism. 

Again,  we may present definability as a form of computability via 
the strangely neglected models from generalised recursion theory. 
For those wanting to dig out the hidden treasures of this 
beautiful and suddenly relevant subject, probably the best 
guide is the 1990 Springer {\em Perspectives in Logic} book of 
Gerald E. Sacks \cite{Sacks:1990} on {\em Higher Recursion Theory}. 

The physical relevance of this extended 
model of computable causality fits nicely with the 
current return to basic notions of key figures concerned with 
quantum gravity and foundational questions in physics. 
 As Lee Smolin  affirms in his book \cite{smolin2006} on {\em The Trouble 
 with Physics} (p.241): ``\dots causality itself is fundamental". 
 
 Smolin is not, of course, thinking of computationally constrained 
 models of causality in the sense of Turing. But there is a convergence 
 of aims and value put on embodiment of higher order relations on 
 very basic structures. Smolin alludes to the work of `early champions' 
 of the role of causality, such as Roger Penrose, 
 Rafael Sorkin, Fay Dowker, and Fotini Markopoulou, and 
 sets out a version of Penrose's ~strong determinism (p.241 again):
 \begin{quote}
 It is not only the case that the spacetime geometry determines what the 
 causal relations are. This can be turned around: Causal relations can 
 determine the spacetime geometry \dots 
ItÕs easy to talk about space or spacetime emerging from something 
 more fundamental, but those who have tried to develop the idea 
 have found it difficult to realize in practice. ... We now believe they 
 failed because they ignored the role that causality plays in 
 spacetime. These days, many of us working on quantum gravity 
 believe that causality itself is fundamental - and is thus meaningful 
 even at a level where the notion of space has disappeared.
 \end{quote}

\section{Embodiment Restored}

The difference between the extended Turing model of computation and 
what is commonly seen from the modellers of process 
  is that there is a proper balance between 
process and information. The embodiment was a key problem for the 
early development of the computer, insufficiently recognised since the 
early days by the 
theorists, fixated on the universality paradigm. 

Rodney Brooks \cite{Brooksta} tells how embodiment 
in the form of stored information 
has re-emerged in AI:
\begin{quote}
Modern researchers are now seriously investigating the embodied 
approach to intelligence and have rediscovered the importance of 
interaction with people as the basis for intelligence. My own work for the 
last twenty five years has been based on these two ideas.
\end{quote}

The mathematics of the extended Turing model is notoriously 
for its technicality. And the mathematical character of the 
global structure based on it is disappointingly pathological 
from the point of view of mathematical 
expectations. But if we expect to model the emergence of global 
relations in terms of local structure -- say that of large-scale structures 
in the real universe, or even of relations expressing globally 
observed natural laws -- then the pathology provides the raw material 
for an embodied language within which to talk about 
the complexity of physical forms we discover about us. 
In retrospect, Hartley Rogers was remarkably prescient when 
he asked about the character of the Turing invariant relations 
in a talk entitled {\em Some problems of definability in recursive function theory}, 
at the Tenth Logic Colloquium, at the University of Leicester in 1965. 
Though at that time people had no appreciation of the physical 
significance of the question. It appearing increasingly likely that 
these relations are the  key to pinning down how basic laws and entities 
emerge as mathematical constraints on causal structure. 
One should be aware though of the schematic nature of the 
understanding they may provide. They are the equivalent of maps of the 
landscape provided by satellite scans, while the real 
work of exploring the  substance of down below is at best 
guided and enlightened by an awareness of the overall 
structuring. 

The character of the all-important automorphism group of the 
Turing universe is still to be pinned down. In the way of progress 
towards explanations of such scientific mysteries as the dichotomy between 
classical and quantum reality, and the removal of the need for 
speculative assumptions of `many-worlds' and multiverses, is the 
so-called {\em Bi-interpretability Conjecture}, which has haunted us for nearly 
thirty years. Very roughly speaking, the conjecture says that the 
Turing definable relations are exactly those with information 
content describable in second-order arithmetic. One consequence of a 
positive solution to this problem would be the ruling out of 
non-trivial automorphisms of the Turing universe. And the 
breaking of the current match between the computability-theoretic 
structures and the physical reality we observe. As things are, we see 
and important partial validation of the conjecture, yielding 
a rigid substructure of the Turing universe replete with uniquely 
defined entities -- much like our classical reality we live in: and 
a wildly ill-defined context, with automorphic displacements corresponding 
to pre-measurement quantum ambiguity. Currently, the 
main clues to the outcome of the full Bi-interpretability 
Conjecture are a lack of more than incremental progress towards 
a positive answer for the last fifteen years, and an   
outline of a strategy for constructing a non-trivial  automorphism 
in the public domain since the late 1990s. 

We summarise some features of the mathematics, and refer the reader 
to sources such as \cite{CO03} and \cite{Co09a} for further detail:
\begin{itemize}
\item Embodiment invalidating the `machine as data' and universality paradigm. 
\item The organic linking of mechanics and emergent outcomes delivering a 
clearer model of supervenience of mentality on brain functionality, and a 
reconciliation of different levels of effectivity. 
\item A reaffirmation of the importance of experiment and evolving hardware, 
for both AI and extended computing generally. 
\item The validating of a route to {\it creation} of new information through 
interaction and emergence.
\item The significance  of definability and its breakdown for the physical 
universe giving a route to the determination of previously unexplained 
constants in the standard model of physics, and of quantum ambiguity 
in terms of a breakdown of definability, so making the multiverse 
redundant as an explanatory tool, and \dots 
\item \dots work by Slaman and Woodin and others 
on establishing partial rigidity of the Turing universe 
(see \cite{AFta}) 
promising an explanation of the existence of our `quasi-classical' 
universe. 

\end{itemize}

As for building intelligent machines, we give the last word to 
Danny Hillis, quoted by Mark Williams in {\it Red Herring} magazine in  
April 03, 2001:
\begin{quote}
I used to think we'd do it by engineering. Now I believe we'll 
evolve them. We're likely to make thinking machines before we 
understand how the mind works, which is kind of backwards.
\end{quote}
So, Turing's computational modelling is showing good signs 
of durability --  but within a turbulent natural environment, 
embodied, full of emergent wonders, and exhibiting a 
computational structure reaching far into unknown regions 
 -- via a revived type structure, the relevance 
 of which is currently neglected by both mathematicians,  
  computer scientists, and beyond. But, as I write this at the 
  start of the centenary year of Alan Turing's birth, there 
  is every sign of a revival of the basic approach 
  which underlay the great discoveries of the first half of the 
  twentieth century. 

Turing, as we know, anticipated much 
of what we are still clarifying about how the world computes.  
And, we hope, he is smiling down on us at the 
centenary of his birth, in a Little Venice nursing home (now the 
Colonnade Hotel).

\bibliographystyle{abbrv} 
\bibliography{challenges.bib}  

\begin{thebibliography}{10}

\bibitem{AFta}
K.~Ambos-Spies and P.~Fejer.
\newblock Degrees of unsolvability.
\newblock {\em unpublished}, 2006.

\bibitem{BBY07}
I.~Binder, M.~Braverman, and M.~Yampolsky.
\newblock Filled {J}ulia sets with empty interior are computable.
\newblock {\em Foundations of Computational Mathematics}, 7(4):405--416, 2007.

\bibitem{Brooks01}
R.~Brooks.
\newblock The relationship between matter and life.
\newblock {\em Nature}, 409:409--411, 2001.

\bibitem{Brooksta}
R.~Brooks.
\newblock The case for embodied intelligence.
\newblock In S.~B. Cooper and J.~van Leeuwen, editors, {\em Alan Turing - His
  Work and Impact}. Elsevier Science, to appear.

\bibitem{CS08}
C.~S. Calude and K.~Svozil.
\newblock Quantum randomness and value indefiniteness.
\newblock {\em Advanced Science Letters}, 1:165--168, 2008.

\bibitem{Chaitinip}
G.~J. Chaitin.
\newblock Metaphysics, metamathematics and metabiology.
\newblock In H.~Zenil, editor, {\em Randomness Through Computation: Some
  Answers, More Questions}. World Scientific, Singapore, 2011.

\bibitem{Co09a}
S.~B. Cooper.
\newblock Emergence as a computability-theoretic phenomenon.
\newblock {\em Applied Mathematics and Computation}, 215(4):1351--1360, 2009.

\bibitem{CO03}
S.~B. Cooper and P.~Odifreddi.
\newblock Incomputability in nature.
\newblock In S.~B. Cooper and S.~S. Goncharov, editors, {\em Computability and
  Models: Perspectives East and West}, pages 137--160. Plenum, New York, 2003.

\bibitem{Damasio1999}
A.~R. Damasio.
\newblock {\em The Feeling of What Happens: Body and Emotion in the Making of
  Consciousness}.
\newblock Harcourt Brace and Co, 1999.

\bibitem{Denning:2010}
P.~J. Denning.
\newblock Ubiquity symposium `{What} is computation?': Opening statement.
\newblock {\em Ubiquity}, 2010.

\bibitem{Deutsch85}
D.~Deutsch.
\newblock Quantum theory, the {C}hurch-{T}uring principle, and the universal
  quantum computer.
\newblock {\em Proc. Royal Soc.}, A400:97--117, 1985.

\bibitem{Fortnow:2010}
L.~Fortnow.
\newblock Ubiquity symposium `{What }is computation?': The enduring legacy of
  the turing machine.
\newblock {\em Ubiquity}, 2010, Dec. 2010.

\bibitem{Kim2000}
J.~Kim.
\newblock {\em Mind in a Physical World: An Essay on the Mind-Body Problem and
  Mental Causation}.
\newblock MIT Press, 2000.

\bibitem{Kim2005}
J.~Kim.
\newblock {\em Physicalism, or Something Near Enough}.
\newblock Princeton University Press, 2005.

\bibitem{Murray02}
J.~D. Murray.
\newblock {\em Mathematical Biology: I. An Introduction}.
\newblock Springer, New York, 3rd edition, 2002.

\bibitem{Penrose2002}
R.~Penrose.
\newblock {\em The Emperor's New Mind: Concerning Computers, Minds, and the
  Laws of Physics (Popular Science)}.
\newblock Oxford University Press, USA, 2002.

\bibitem{pinker}
S.~Pinker.
\newblock {\em How the Mind Works}.
\newblock W.W. Norton, New York, 1997.

\bibitem{Put60}
H.~Putnam.
\newblock Minds and machines (1960).
\newblock In Putnam, editor, {\em Mind, Language, and Reality}, pages 362--385.
  Cambridge University Press, 1975.

\bibitem{Sacks:1990}
G.~E. Sacks.
\newblock {\em Higher recursion theory}.
\newblock Springer-Verlag New York, Inc., New York, NY, USA, 1990.

\bibitem{Slo09}
A.~Sloman.
\newblock Some requirements for human-like robots: Why the recent over-emphasis
  on embodiment has held up progress.
\newblock In B.~Sendhoff et~al., editors, {\em Creating Brain-Like
  Intelligence}, volume 5436 of {\em LNCS}, pages 248--277. Springer, 2009.

\bibitem{Smolensky88}
P.~Smolensky.
\newblock On the proper treatment of connectionism.
\newblock {\em Behavioral and Brain Sciences}, 11:1--74, 1988.

\bibitem{smolin2006}
L.~Smolin.
\newblock {\em The trouble with physics: the rise of string theory, the fall of
  a science, and what comes next}.
\newblock Houghton Mifflin Co., 2006.

\bibitem{Ta07}
N.~N. Taleb.
\newblock {\em The Black Swan: The Impact of the Highly Improbable}.
\newblock Random House, 2010.

\bibitem{Teuscher02}
C.~Teuscher.
\newblock {\em Turing's Connectionism. An Investigation of Neural Network
  Architectures}.
\newblock Springer-Verlag, London, 2002.

\bibitem{Turing1939}
A.~M. Turing.
\newblock Systems of logic based on ordinals.
\newblock {\em Proceedings of the London Mathematical Society}, 45:161--228,
  1939.
\newblock reprinted in Alan M. Turing, Collected Works: Mathematical Logic, pp.
  81-148.

\bibitem{Turing1950}
A.~M. Turing.
\newblock Computing machinery and intelligence.
\newblock {\em MIND}, 59(236):433--460, Oct. 1950.

\bibitem{Tu52}
A.~M. Turing.
\newblock The chemical basis of morphogenesis.
\newblock {\em Phil. Trans. of the Royal Society of London. Series B,
  Biological Sciences}, 237(641):37--72, 1952.

\bibitem{WG2005}
P.~Wegner and D.~Goldin.
\newblock The {C}hurch-{T}uring {T}hesis: Breaking the myth.
\newblock In S.~B. Cooper and B.~L\"{o}we, editors, {\em CiE 2005: New
  Computational Paradigms}, volume 3526 of {\em LNCS}. Springer, 2005.

\end{thebibliography}

\end{document}